\documentclass[12pt,psamsfonts]{amsart}
\usepackage{amsmath,amsthm,amsfonts,amssymb}

\usepackage{eucal}
\newcommand{\1}{\mathbf{1}}

\newcommand{\M}{\mathcal{M}}
\newcommand{\K}{\mathcal{K}}

\newcommand{\lan}{\langle}
\newcommand{\ra}{\rangle}

\newcommand{\T}{\mathcal{T}}

\newcommand{\F}{\mathcal{F}}

\newcommand{\ga}{\gamma}

\newcommand{\g}{\mathfrak{g}}

\newcommand{\Q}{Q}

\newcommand{\Cg}{\CC(\g,q,\ell)}
\newcommand{\la}{{\lambda}}

\DeclareMathOperator{\U}{U}
\DeclareMathOperator{\rank}{rank} 
\DeclareMathOperator{\End}{End} \DeclareMathOperator{\GL}{GL}
\DeclareMathOperator{\Hom}{Hom}

\DeclareMathOperator{\SL}{SL}\DeclareMathOperator{\SU}{SU}\DeclareMathOperator{\PSL}{PSL}

\newcommand{\C}{\mathbb C}
\newcommand{\CC}{\mathcal{C}}
\newcommand{\CE}{\mathcal{E}}
\newcommand{\CH}{\mathcal{H}}
\newcommand{\mZ}{\mathcal{Z}}
\newcommand{\Z}{\mathbb Z}
\newcommand{\ot}{\otimes}
\newcommand{\B}{\mathcal{B}}

\numberwithin{equation}{section}

\newtheorem{theorem}[equation]{Theorem}
\newtheorem{defn}[equation]{Definition}

\newtheorem{lemma}[equation]{Lemma}

\newtheorem{prop}[equation]{Proposition}
\theoremstyle{definition}

\newtheorem{remark}[equation]{Remark}

\begin{document}
\title[Exotic modular categories]
{On exotic modular tensor categories}

\author{Seung-moon Hong}
\email{seuhong@indiana.edu}
\address{Department of Mathematics\\
    Indiana University \\
    Bloomington, IN 47405\\
    U.S.A.}

\author{Eric Rowell}
\email{rowell@math.tamu.edu}
\address{Department of Mathematics\\
    Texas A\&M University \\
    College Station, TX 77843\\
    U.S.A.}

\author{Zhenghan Wang}
\email{zhenghwa@microsoft.com}
\address{Microsoft Station Q\\Elings Hall 2237\\
    University of California\\
    Santa Barbara, CA 93106\\
    U.S.A.}

\dedicatory{Dedicated to the memory of Xiao-Song Lin}

\thanks{The first and third authors are partially supported by NSF FRG grant DMS-034772.
The third author likes to thank F.~Xu, M.~ M\"uger, V. Ostrik, and
Y.-Z.~ Huang for helpful correspondence.}

\begin{abstract}

It has been conjectured that every $(2+1)$-TQFT is a
Chern-Simons-Witten (CSW) theory labeled by a pair $(G,\lambda)$,
where $G$ is a compact Lie group, and $\lambda \in H^4(BG;\Z)$ a
cohomology class.  We study two TQFTs constructed from Jones'
subfactor theory which are believed to be counterexamples to this
conjecture: one is the quantum double of the even sectors of the
$E_6$ subfactor, and the other is the quantum double of the even
sectors of the Haagerup subfactor. We cannot prove mathematically
that the two TQFTs are indeed counterexamples because CSW TQFTs,
while physically defined, are not yet mathematically constructed
for every pair $(G,\lambda)$.  The cases that are constructed
mathematically include:

\begin{enumerate}

\item  $G$ is a finite group---the Dijkgraaf-Witten TQFTs;

\item  $G$ is torus $T^n$;

\item $G$ is a connected semi-simple Lie group---the
Reshetikhin-Turaev TQFTs.

\end{enumerate}

We prove that the two TQFTs are not among those mathematically
constructed TQFTs or their direct products.  Both TQFTs are of the
Turaev-Viro type: quantum doubles of spherical tensor categories.
We further prove that neither TQFT is a quantum double of a
braided fusion category, and give evidence that neither is an
orbifold or coset of TQFTs above. Moreover, representation of the
braid groups from the half $E_6$ TQFT can be used to build
universal topological quantum computers, and the same is expected
for the Haagerup case.

\end{abstract}
\maketitle

\section{Introduction}
\label{s:intro}

In his seminal paper \cite{Witten89}, E.~Witten invented
Chern-Simons $(2+1)$-topological quantum field theory (TQFT), and
discovered a relation between Chern-Simons TQFTs and
Wess-Zumino-Novikov-Witten (WZW) conformal field theories (CFTs).
To be more precise, CFTs here should be referred to as chiral
CFTs, as opposed to full CFTs. The connection between
Chern-Simons-Witten $(2+1)$-TQFTs and WZW models has spawned an
application of TQFT and rational CFT (RCFT) to condensed matter
physics (see \cite{RSW} and the references therein). In fractional
quantum Hall liquids, Chern-Simons-Witten theories are used to
describe emerged topological properties of the bulk electron
liquids, whereas the corresponding CFTs describe the boundary
physics of the Hall liquids (see \cite{Wilczek} and the references
therein). A unifying theme in the mathematical formulation of both
$(2+1)$-TQFTs and CFTs is the notion of a modular tensor category
(MTC) \cite{Turaev}. Modular tensor categories are the algebraic
data that faithfully encode $(2+1)$-TQFTs \cite{Turaev}, and are
used to describe anyonic properties of certain quantum systems
(see \cite{Kitaev2} \cite{DFN} \cite{Wang} and the references
therein). In this paper, we will use the terms $(2+1)$-TQFT, or
just TQFT in the future, and MTC interchangeably (We warn readers that it is
an open question whether or not TQFTs and MTCs are in one-one correspondence,
see e.g. \cite{BK}.
But an MTC gives rise to a unique TQFT \cite{Turaev}).
Our interest in MTCs comes from topological quantum computing by braiding
non-abelian anyons in the sense of \cite{FKLW}
(cf.\cite{Kitaev1}). From this perspective, we are interested in
an abstract approach to MTCs free of algebraic structures such as
vertex operator algebras (VOAs) or local conformal nets of von
Neumann algebras, whose representation theory gives rise to MTCs
(see \cite{Hu}\cite{KLM}\cite{Evans} and the references therein).

Known examples of MTCs that are realized by anyonic quantum
systems in real materials are certain abelian MTCs encoding
Witten's quantum Chern-Simons theories for abelian gauge groups at
low levels (see \cite{Wilczek}). The physical systems are
$2$-dimensional electron liquids immersed in strong perpendicular
magnetic fields that exhibit the so-called fractional quantum Hall
effect (FQHE). In these physical systems, the representations of
the braid groups from the MTCs describe braiding statistics of the
quasi-particles, which are neither bosons nor fermions. F.~Wilczek
named such exotic quasi-particles \textit{anyons}. Confirmation of
the realization of non-abelian MTCs in FQH liquids is pursued
actively in experiments (see \cite{DFN} and the references
therein).

Inspired by FQHE, we may imagine that there are physical systems
to realize many MTCs. With this possibility in mind, we are
interested in the construction and classification of MTCs.  Since
TQFTs and CFTs are closely related to each other, we may expect
all the known constructions of new CFTs from given CFTs such as
coset, orbifold, and simple current extension can be translated
into the TQFT side, and then to the MTC side in a purely
categorical way. After many beautiful works, it seems that those
constructions cannot in general be defined in the purely
categorical setting. On the CFT side, it has been expressed
several times in the literature that all known rational CFTs are
covered by a single construction: Witten's quantum Chern-Simons
theory.  In particular, the following conjecture is stated in
\cite{MS1}:

\vspace{0.2in}

{\bf Conjecture 1:} {\it The modular functor of any unitary RCFT
is equivalent to the modular functor of some Chern-Simons-Witten
{\bf (CSW)} theory defined by the pair $(G,\lambda)$ with $G$ a
compact group and $\lambda\in H^4(BG;\Z).$ }

\vspace{0.2in}

Another conjecture, attributed to E.~Witten \cite{Witten3}, was
stated as Conjecture 3 in \cite{MS1}:

\vspace{0.2in}

{\bf Conjecture 3:} {\it All three dimensional topological field
theories are CSW theory for some appropriate (super)-group.}

\vspace{0.2in}

CSW theories with compact Lie groups are written down in
\cite{DijkgraafWitten}, and they are labeled by a pair
$(G,\lambda)$, where $G$ is a compact Lie group, and $\lambda\in
H^4(BG;\Z)$.  A modular functor is just the $2$-dimensional part
of a TQFT \cite{Turaev}\cite{MS2}; TQFTs from CSW theory as in the
conjectures will be called {\bf CSW} TQFTs. Therefore, we
paraphrase the two conjectures as:

\vspace{0.2in}

{\bf Conjecture CSW:} {\it Every $(2+1)$-TQFT is a CSW TQFT for
some pair $(G,\lambda)$, where where $G$ is a compact Lie group,
and $\lambda\in H^4(BG;\Z)$ a cohomology class.}

\vspace{0.2in}

Since TQFTs are faithfully encoded by MTCs \cite{Turaev},
translated into the MTC side, this conjecture says that any
unitary MTC is equivalent to one from some unitary CSW TQFT. If
this conjecture holds, we will have a conceptual classification of
$(2+1)$-TQFTs. Of course even if the conjecture were true, to make
such a classification into a mathematical theorem is still very
difficult.

There are three families of compact Lie groups for which we have
mathematical realizations of the corresponding $(2+1)$ CSW TQFTs:

\begin{enumerate}

\item $G$ is finite \cite{DijkgraafWitten}\cite{FreedQuinn};

\item $G$ is a torus $T^n$ \cite{Manoliu}\cite{BelovMoore};

\item $G$ is a connected semi-simple Lie group
\cite{ReshetikhinT}\cite{Turaev}.

\end{enumerate}

Given such an attractive picture, we are interested in the
question whether or not all known TQFTs fit into this framework.
An MTC or TQFT will be called {\it exotic} if it cannot be
constructed from a CSW theory. In this paper, we will study two
MTCs which seem to be exotic: the quantum doubles $\mZ(\CE)$ and
$\mZ(\CH)$ of the spherical fusion categories $\CE$ and $\CH$
generated by the even sectors of the $E_6$ subfactor (a.k.a.
$\frac{1}{2}E_6$), and the even sectors of the Haagerup subfactor
of index $\frac{5+\sqrt{13}}{2}$ \cite{AsaedaH}. Unfortunately we
cannot prove that these two unitary MTCs are indeed exotic. The
difficulty lies in describing mathematically all unitary CSW MTCs,
in particular those from non-connected, non-simply-connected Lie
groups $G$. When $G$ is finite, the corresponding Dijkgraaf-Witten
TQFTs are (twisted) quantum doubles of the group categories, which
are well understood mathematically (see \cite{BK}).  When $G$ is a
torus, the corresponding TQFTs are abelian, and are classified in
\cite{BelovMoore}. When $G$ is a connected semi-simple Lie group,
the CSW MTCs are believed to correspond mathematically to MTCs
constructed by N.~Reshetikhin and V.~Turaev based on the
representation theory of quantum groups
\cite{ReshetikhinT}\cite{Turaev}.  We will see that the two
seemingly exotic MTCs cannot be constructed by using $G$ finite or
$G$ a torus. Therefore, we will study whether or not they can be
obtained from categories constructed from quantum groups.

Quantum groups are deformations of semi-simple Lie algebras.  The
standard quantum group theory does not have a well established
theory to cover non-connected Lie groups.  So our translation of
Conjecture CSW to quantum group setting is not faithful since we
will only study MTCs from deforming semi-simple Lie algebras. MTCs
constructed in this way will be called {\bf quantum group
categories} in this paper, which are constructed mathematically
(see \cite{Turaev} \cite{BK} and the references therein).  To
remedy the situation to some extent, we will consider coset and
orbifold constructions from quantum group categories in Section
\ref{coset}. It is known that coset and orbifold theories are
included in the CSW theories by using appropriate compact Lie
groups, in particular non-connected Lie groups \cite{MS1}.

There are new methods to construct MTCs. In particular many
examples are constructed through VOAs and von Neumann algebras.
Several experts in the mathematical community believe that those
examples contain new MTCs that are not CSW MTCs. But as alluded
above to prove such a statement is mathematically difficult. First
even restricted to quantum group categories, the mathematical
characterization of all MTCs from quantum group categories plus
coset, orbifold, and simple current extension is hard, if not
impossible. Secondly, the potentially new examples of MTCs are
complicated measured by the number of simple object types. Another
construction of MTCs is the quantum double, which is a categorical
generalization of the Drinfeld double of quantum groups.  Such
MTCs give rise to TQFTs of the Turaev-Viro type \cite{TV} and
naturally arise in subfactor theory by A.~Ocneanu's asymptotic
inclusions construction (see \cite{Evans}). The categorical
formulation of Ocneanu's construction is M.~M\"uger's beautiful
theorem that the quantum double of any spherical category is an
MTC \cite{MugerII}. The authors do not know how to construct
quantum double TQFTs from CSW theory in general, except for the
finite group case; hence general quantum double TQFTs are
potentially exotic, and might be the only exotic ones.  Maybe
quantum double TQFTs can be constructed as CSW theory for some
appropriate super-groups as Witten conjectured, but we are not
aware of such mathematical theories.

The most famous examples of double TQFTs are related to the
Haagerup subfactor of index $\frac{5+\sqrt{13}}{2}$. The Ocneanu
construction gives rise to a unitary MTC which is the quantum
double of a spherical category of $10$ simple object types. This
spherical category of $10$ simple object types is not braided
because there are simple objects $X,Y$ such that $X\otimes Y$ is
not isomorphic to $Y\otimes X$. Since the Haagerup subfactor
cannot be constructed from quantum groups \cite{AsaedaH}, it is
unlikely that the corresponding MTC is isomorphic to an MTC from
quantum groups or a quantum double of a quantum group category.
Recall that the standard construction of quantum group categories
are always braided. But this MTC is difficult to study explicitly
as the number of the simple object types, called the {\it rank} of
an MTC, is big. An alternative is to study the double of the even
sectors of the Haagerup subfactor. The even sectors form a
non-braided spherical category with $6$ simple types; its quantum
double is a unitary MTC of rank $12$, which will be called {\bf
the Haagerup MTC}. There is another simpler category which has
similar exoticness: the quantum double of $\frac{1}{2}E_6$, which
is of rank $10$. But note that a quantum double of a non-braided
spherical category can be constructed by CSW theory sometimes. For
example, if we double the group category $S_3$ of rank $6$, which
is not braided, we get a unitary MTC of rank $8$, which is a
Dijkgraaf-Witten MTC.

The spherical category $\frac{1}{2}E_6$ was brought to our
attention by V.~Ostrik \cite{Ostrik1}, and its double is worked
out in \cite{BEY} \cite{Iz}. In V.~Ostrik's paper on the
classification of rank=$3$ fusion categories with braidings, he
conjectured that there is only one set of fusion rules of rank=$3$
without braidings. This set of fusion rules is realized by the
$\frac{1}{2}E_6$ fusion rules, and is known to be non-braided. If
we denote by $1,x,y$ three representatives of the simple object
types, their fusion rules are: $x^2=1+2x+y,xy=yx=x,y^2=1$.  Note
that we simply write tensor product as multiplication and will
denote the $\frac{1}{2}E_6$ category as $\CE$, and its double
$\mZ(\CE)$.  If we denote by
$1,\alpha,\alpha^*,\rho,{}_{\alpha}\rho,{}_{\alpha^*}\rho$ six
representatives of the even sectors of the Haagerup subfactor,
their fusion rules are $\alpha \alpha^*=1, \alpha^2=\alpha^*,
\alpha\rho={}_{\alpha}\rho, \alpha^*\rho={}_{\alpha^*}\rho,
\alpha\rho=\rho\alpha^*,\rho^2=1+\rho+{}_{\alpha}\rho+{}_{\alpha^*}\rho.$
We will denote this rank $6$ unitary fusion category by $\CH$, and
its double $\mZ(\CH)$, the Haagerup MTC.

Our main result is:

\begin{theorem}
\label{Main Thm}

Let $\CE, \CH$ be the non-braided unitary spherical categories
above, and $\mZ(\CE), \mZ(\CH)$ be their quantum doubled MTCs.
Then $\mZ(\CE), \mZ(\CH)$

\begin{enumerate}

\item are prime, i.e. there are no non-trivial modular
subcategories; hence are not a product of two MTCs;

\item have non-integral global quantum dimension $D^2$, hence are
not CSW MTCs for finite or torus Lie group $G$;

\item are not quantum group categories;

\item are not quantum doubles of any braided fusion categories;

\item give rise to representations of $\SL(2,\Z)$ that factor over
a finite group.

\item $\mZ(\CE)$ gives rise to representations of the braid groups
with infinite images;

\item $\mZ(\CE)$ is a decomposable bimodule category over a
pre-modular subcategory of quantum group type, but $\mZ(\CH)$ has no
such decompositions.

\end{enumerate}

\end{theorem}

It is known that all quantum double MTCs are non-chiral in the
sense their topological central charges are $0$, hence anomaly
free in the sense that the representations afforded with all
mapping class groups are linear representations rather than
projective ones. This is a subtle point since the topological
central charge is $0$ does not imply the chiral central charge of
the corresponding CFT, if there is one, is $0$.  The topological
central charge is only defined modulo $8$, so topological central
charge being $0$ means the chiral central charge of the
corresponding CFT, if it exists, is $0$ mod $8$.
 It is possible that $\mZ(\CE)$ or $\mZ(\CH)$ can be constructed as
an orbifold of a quantum group category with topological central
charge $=0$ mod $8$ or as a coset category of quantum group
category.  But as we will see in Section~\ref{coset}, an MTC of an
orbifold CFT has global quantum dimension at least $4$ times of
that of the the original MTC, so it is unlikely for either to be
an orbifold.  The coset construction is more complicated, but the
constraint of central charge being multiples of $8$ restricts the
possible cosets significantly.  We will leave a detailed analysis
to the future.

To prove that any of the two MTCs are indeed not CSW MTCs for some
pair $(G,\lambda)$, we need a classification of all CSW CFTs of
central charge $0$ mod $8$ and $10$ or $12$ primary fields.  If
the classification is simple enough, we may just examine the list
to show that our exotic examples are not among the associated
MTCs. This seems to be difficult.

The existence of a pre-modular subcategory in $\mZ(\CE)$ raises an
interesting possibility.  The tensor sub-category generated by
$X_4$ has $6$ simple objects: $$\{\mathbf{1},Y,X_4,X_5,U,V\}.$$
The Bratteli diagram for decomposing tensor powers of $X_4$ is
identical to that of a pre-modular category associated with the
subcategory of non-spin representations of quantum
$\mathfrak{so}_3$ at a $12$th root of unity (see
Prop.~\ref{bimodule}). This suggests the possibility that
$\mZ(\CE)$ might be related to an $O(3)$-CSW MTC.  Recall that
$H^4(BO(3);\Z)=\Z\oplus \Z_2$, so for each level $k$, there are
two CSW MTCs.  It will be interesting to compare the $12$th root
of unity $O(3)$-CSW theory for the nontrivial $\Z_2$ level with
$\mZ(\CE)$. A construction in CFT that we will not consider is the
simple current extension.  It is possible that $\mZ(\CE)$ is a
simple current extension of a quantum group category, or a
sub-category of the simple extension of a quantum group category.
On the other hand, $\mZ(\CH)$ has no sub tensor categories, so the
above discussion seems not applicable to $\mZ(\CH)$.  It is still
possible that $\mZ(\CE)$ or $\mZ(\CH)$ can be constructed as the
quantum double of a spherical quantum group category which is not
braided. As far as we know there are no systematic ways to produce
spherical quantum group categories that are not ribbon.  But every
fusion category comes from a \emph{weak Hopf algebra}, therefore there will be
no exotic MTCs if the term \emph{quantum group} is too liberal
\cite{Ostrik2}.

A sequence of potentially new chiral CFTs were recently
constructed in \cite{KL}\cite{Xu1}. The associated MTCs might be
exotic. The simplest one in this sequence has an MTC equivalent to
the mirror MTC of $SU(5)_1\times SO(7)_1$ with chiral central
charge $16.5$. It will be interesting to analyze these categories.

The paper is organized as follows.  In Section \ref{s:double}, we
calculate the $S,T$ matrices, and fusion rules of both TQFTs. Then
we deduce several observations including (7) of Theorem \ref{Main
Thm}.  In Section \ref{s:appl}, we prove (3) of Theorem \ref{Main
Thm} and neither $\mZ(\CE)$ nor $\mZ(\CH)$ is a product of two
MTCs. Then both MTCs are prime because if there were a nontrivial
modular subcategory of $\mZ(\CE)$ or $\mZ(\CH)$,
\cite{MugerLMS}[Theorem 4.2] implies that $\mZ(\CE)$ would be a
non-trivial product of two MTCs, a contradiction. Sections
\ref{double}, \ref{sl2z}, \ref{braidreps} are devoted to the
proofs of (4) (5) (6) of Theorem \ref{Main Thm}. In Section
\ref{coset}, we give evidence that neither theory is an orbifold
or coset.  In the appendix, we give an explicit description of the
category $\mZ(\CE)$ from the definition of half braidings.

As a final remark, regardless of the relevance to the Conjecture
CSW, our work seems to be the most detailed study of non-quantum
group TQFTs besides the Dijkgraaf-Witten, abelian and their direct
product TQFTs.  We also understand that some of the results are
well-known to some experts, but are not well-documented.

\section{Categories $\mZ(\CE)$ and $\mZ(\CH)$}
\label{s:double}

\subsection{$\mZ(\CE)$}

The spherical category $\CE$ is studied in \cite{BEY} \cite{Iz},
and the associated Turaev-Viro invariant is studied in \cite{SW}.
All spherical categories with the same set of fusion rules are
worked out in great detail in \cite{HH}. Those categories have
three isomorphism classes of simple objects denoted by $1,x,y$,
their fusion rules are: $x^2=1+2x+y,xy=yx=x,y^2=1$. The categories
are called $\frac{1}{2}E_6$ because the fusion rules can be
encoded by half of the Dynkin diagram $E_6$. There is an
essentially unique unitary spherical category with this set of
fusion rules up to complex conjugation. We pick the same one as in
\cite{HH} as our $\CE$, and all conceptual conclusions will be
same for other choices except when specific complex parameters are
involved. By direct computation from the definition (details are
given in the appendix), we find that the quantum double of $\CE$
has $10$ simple object types, of which $8$ are self-dual and the
other two are dual to each other. We label the $10$ simple objects
by $\mathbf{1}:=(1,e_1)$, $Y:=(y,e_y)$, $X_i:=(x,e_{x_i})$ for
$i=1,2,\cdots,5$, $U:=(1+x,e_{1+x})$, $V:=(y+x,e_{y+x})$, and
$W:=(1+y+x,e_{1+y+x})$, where the half-braiding notations as in
the appendix are used here.

\subsubsection{S-matrix, T-matrix}\label{SandT}

Once we have the list of all simple objects as half braidings as
in the appendix, it is easy to compute the $S$-matrix, and the
$T$-matrix. (The $S,T$ matrices can also be computed from
\cite{Iz} and is also contained in \cite{EP1}.)  In an MTC,
we have $D=\sqrt{\sum_i d_i^2}$, where $i$ goes over all simple
object types.  There are various names in the literature for $D$ and $D^2$.  We call $D^2$ the
global quantum dimension, and $D$ the total quantum order.  Total quantum order is
not a standard terminology,
and is inspired by the role that $D$ plays in topological entropy for topological
phases of matter.  The $S$-matrix is $S=\frac{1}{D}\tilde{S}$, where the
total quantum order $D=\sqrt{\dim(\mZ(\CE))}=6+2\sqrt{3}$, and
$\tilde{S}$ is as follows:

$\left(\begin{smallmatrix}1&1&\sqrt{3}+1&\sqrt{3}+1&\sqrt{3}+1&\sqrt{3}+1&\sqrt{3}+1&\sqrt{3}+2&\sqrt{3}+2&\sqrt{3}+3\\
1&1&-\sqrt{3}-1&-\sqrt{3}-1&-\sqrt{3}-1&\sqrt{3}+1&\sqrt{3}+1&\sqrt{3}+2&\sqrt{3}+2&-\sqrt{3}-3\\
\sqrt{3}+1&-\sqrt{3}-1&0&0&0&2(\sqrt{3}+1)&-2(\sqrt{3}+1)&-\sqrt{3}-1&\sqrt{3}+1&0\\
\sqrt{3}+1&-\sqrt{3}-1&0&-i(\sqrt{3}+3)&i(\sqrt{3}+3)&-\sqrt{3}-1&\sqrt{3}+1&-\sqrt{3}-1&\sqrt{3}+1&0\\
\sqrt{3}+1&-\sqrt{3}-1&0&i(\sqrt{3}+3)&-i(\sqrt{3}+3)&-\sqrt{3}-1&\sqrt{3}+1&-\sqrt{3}-1&\sqrt{3}+1&0\\
\sqrt{3}+1&\sqrt{3}+1&2(\sqrt{3}+1)&-\sqrt{3}-1&-\sqrt{3}-1&\sqrt{3}+1&\sqrt{3}+1&-\sqrt{3}-1&-\sqrt{3}-1&0\\
\sqrt{3}+1&\sqrt{3}+1&-2(\sqrt{3}+1)&\sqrt{3}+1&\sqrt{3}+1&\sqrt{3}+1&\sqrt{3}+1&-\sqrt{3}-1&-\sqrt{3}-1&0\\
\sqrt{3}+2&\sqrt{3}+2&-\sqrt{3}-1&-\sqrt{3}-1&-\sqrt{3}-1&-\sqrt{3}-1&-\sqrt{3}-1&1&1&\sqrt{3}+3\\
\sqrt{3}+2&\sqrt{3}+2&\sqrt{3}+1&\sqrt{3}+1&\sqrt{3}+1&-\sqrt{3}-1&-\sqrt{3}-1&1&1&-\sqrt{3}-3\\
\sqrt{3}+3&-\sqrt{3}-3&0&0&0&0&0&\sqrt{3}+3&-\sqrt{3}-3&0
\end{smallmatrix}
\right)$

\vspace{.1in}

The $T$-matrix is diagonal with diagonal entries $\theta_{x_i}, i \in
\Gamma$. The following are the entries:

\vspace{.1in}

$\theta_{\mathbf{1}}=1$, $\theta_{Y}=-1$, $\theta_{X_1}=-i$,
$\theta_{X_2}=\theta_{X_3}=e^{5\pi i/6}$, $\theta_{X_4}=e^{\pi
i/3}$, $\theta_{X_5}=e^{-2\pi i/3}$, $\theta_{U}=1$,
$\theta_{V}=-1$, and $\theta_{W}=1$.

Quantum dimensions of simple objects are among
$\{1,1+\sqrt{3},2+\sqrt{3},3+\sqrt{3}\}$. Twists of simple objects
are all $12^{\rm{th}}$ root of unity.  Notice that every simple
object is self-dual except that $X_2$ and $X_3$ which are dual to
each other.

\subsubsection{Fusion rules}

The fusion rules for $\mZ(\CE)$ can be obtained from the
$S$-matrix via the Verlinde formula or more directly from the
half-braidings. Set $F=\{\1,Y,X_4,X_5,U,V\}$ and
$M=\{X_1,X_2,X_3,W\}$. We record the non-trivial rules in the
following:

{\tiny
$$\begin{array}{|c|c|c|c|c|c|}
\hline F\ot F & Y& X_4 & X_5 & U & V\\
\hline
Y& \1 &X_5 &X_4 & V & U \\
\hline
X_4& X_5& \1+X_4+V & Y+X_5+U & X_5+U+V&X_4+U+V\\
\hline
X_5& X_4& Y+X_5+U &\1+X_4+V& X_4+U+V& X_5+U+V\\
\hline
U& V& X_5+U+V & X_4+U+V& \1+X_4+X_5+U+V& Y+X_4+X_5+U+V\\
\hline
V& U & X_4+U+V&X_5+U+V & Y+X_4+X_5+U+V & \1+X_4+X_5+U+V\\
\hline
\end{array}$$
$$\begin{array}{|c|c|c|c|c|}
\hline F\ot M & X_1 & X_2 & X_3 & W\\
\hline
Y& X_1& X_3 &X_2& W\\
\hline
X_4& X_1+W& X_3+W& X_2+W& X_1+X_2+X_3+W\\
\hline
X_5& X_1+W& X_2+W & X_3+W &X_1+X_2+X_3+W\\
\hline
U& X_2+X_3+W &X_1+X_3+W& X_1+X_2+W & X_1+X_2+X_3+2W\\
\hline
V& X_2+X_3+W& X_1+X_3+W &X_1+X_3+W & X_1+X_2+X_3+2W\\
\hline
\end{array}$$

}

{\tiny
$$
\begin{array}{|c|c|c|c|c|}
\hline M\ot M & X_1 & X_2 & X_3 & W\\
\hline
X_1& \1+Y+X_4+X_5 & U+V &U+V& X_4+X_5+U+V\\
\hline
X_2& U+V& Y+X_4+U&\1+X_5+V& X_4+X_5+U+V\\
\hline
X_3& U+V& \1+X_5+V& Y+X_4+U&X_4+X_5+U+V\\
\hline
W& X_4+X_5+U+V& X_4+X_5+U+V&X_4+X_5+U+V& \1+Y+X_4+X_5+2U+2V \\
\hline
\end{array}
$$
}

From the fusion rules, we observe the following:

\begin{prop}\label{bimodule}

$\mZ(\CE)$ is a decomposable bimodule category over a pre-modular
subcategory of quantum group type.

\end{prop}

\begin{proof}

Observe that the tensor subcategory $\F$ generated by $X_4$ has
$6$ simple objects; namely the simple objects in the set $F$
above. The Bratteli diagram for decomposing tensor powers of $X_4$
is identical to that of a pre-modular category associated with the
subcategory of non-spin representations of quantum
$\mathfrak{so}_3$ at a $12$th root of unity. In fact, the
eigenvalues of the braiding morphism $c_{X_4,X_4}$ are identical
to those of the fusion categories corresponding to BMW-algebras
$BMW_n(q^2,q)$ with $q=e^{\pi i/6}$, and so it follows from the
Tuba-Wenzl classification \cite{TuWe} that these two categories
are braided equivalent.  Moreover, if one takes the semisimple
abelian subcategory $\M$ generated by the simple objects in the
set $M$, one sees that $\mZ(\CE)=\F\oplus\M$ is $\Z_2$-graded with
$\F=\mZ(\CE)_1$ and $\M=\mZ(\CE)_{-1}$ and, moreover, $\M$ is a
bimodule category over $\F$.


\end{proof}

This situation also has interesting connections to Conjecture 5.2
in M\"uger's \cite{MugerLMS}. There he observes (Theorem 3.2) that
if a modular category $\B$ contains a semisimple tensor
sub-category $\K$ then $\dim(\K)\cdot \dim(C_\B(\K))=\dim(\B)$,
where $C_\B(\K)$ is the centralizer subcategory of $\K$ in $\B$.
In the case of $\F\subset \mZ(\CE)$ above, $C_{\mZ(\CE)}(\F)$ is
the subcategory with simple objects $\mathbf{1}$ and $Y$, so that
$\dim(C_{\mZ(\CE)}(\F))=2$. M\"uger calls this a \emph{minimal
modular extension} of $\F$.  He conjectures that any unitary
premodular category $\K$ has a minimal modular extension, that is
$\K\subset \widehat{\K}$ where $\widehat{\K}$ is modular and
$\dim(\widehat{\K})=\dim(\K)\cdot\dim(C_{\widehat{\K}}(\K)).$
Notice also that $\F$ above has at least two such: $\F\subset
\mZ(\CE)$ and $\F\subset\CC(\mathfrak{sl}_2,e^{\pi i/12},12)$.
  This illustrates that the minimal modular extension fails
to be unique, and in fact two such extensions can have different
ranks!

\subsection{$\mZ(\CH)$}
The modular category $\mZ(\CH)$ has rank 12; we label and order
the simple objects as follows:
$\{\1,\pi_1,\pi_2,\sigma_1,\sigma_2,\sigma_3,\mu_1,\cdots,\mu_6\}$
using an abbreviated version of the labeling found in \cite{Iz}.
The quantum dimensions of the non-trivial simple objects are $3d,
3d+1$, and $3d+2$ where $d=\frac{3+\sqrt{13}}{2}$, and the global quantum
dimension is $\dim(\mZ(\CH))=(\frac{39+3\sqrt{13}}{2})^2$.

\subsubsection{S-matrix, T-matrix}\label{HaagSandT}
The modular $S,T$ matrices are also contained in \cite{EP2}.
The review article \cite{E} mentioned
a paper in preparation that contains the explicit expressions of the $x_i$'s below.
The total quantum order
$D=\sqrt{\dim(\mZ(\CH))}=\frac{39+3\sqrt{13}}{2}$, and let $x_i$
denote the six roots of the polynomial
${x}^{6}-{x}^{5}-5\,{x}^{4}+4\,{x}^{3}+6\,{x}^{2}-3\,x-1$ ordered
as follows:

\begin{eqnarray*}&x_1&\approx 0.7092097741,\ x_2\approx 1.497021496,
 \ x_3\approx 1.941883635,\\&x_4&\approx- 0.2410733605,\ x_5\approx- 1.136129493,
\quad \text{and} \quad x_6\approx- 1.770912051.\end{eqnarray*}
With this notation, the $S$-matrix for $\mZ(\CH)$ is
$S=\tilde{S}/D$,
 where:
$$\tilde{S}=\begin{pmatrix} A & B\\
B^T & C
  \end{pmatrix}$$
where $A, B$ and $C$ are the following matrices:
\begin{eqnarray*}
&A&= \begin{pmatrix}1&3d+1&3d+2&3d+2&3d+2&3d+2
\\3d+1&1&3d+2&3d+2&3d+2&3d+2
\\3d+2&3d+2&6d+4&-3d-2&-3d-2&-3d-2
\\3d+2&3d+2&-3d-2&6d+4&-3d-2&-3d-2
\\3d+2&3d+2&-3d-2&-3d-2&-3d-2&6d+4
\\3d+2&3d+2&-3d-2&-3d-2&6d+4&-3d-2
\end{pmatrix},\\
&B&=3d\begin{pmatrix} 1&1&1&1&1&1\\-1&-1&-1&-1&-1&-1
\\0&0&0&0&0&0\\0&0&0&0&0&0
\\0&0&0&0&0&0\\0&0&0&0&0&0
\end{pmatrix},\,
C=3d\begin{pmatrix} x_1&x_3&x_6&x_2&x_4&x_5\\x_3&x_1&
x_2&x_6&x_5&x_4\\
x_6&x_2&x_4&x_5&x_1&x_3
\\x_2&x_6&x_5&x_4&
x_3&x_1\\x_4&x_5&x_1&x_3&x_6&x_2\\x_5&
x_4&x_3&x_1&x_2&x_6\end{pmatrix}
\end{eqnarray*}
Since the entries of $S$ are all real numbers, a simple argument
using the Verlinde formulas shows that all objects in $\mZ(\CH)$
are self-dual.

Now fix $\ga=e^{2\pi i/13}$. The T-matrix is given in \cite{Iz}
and is the diagonal matrix with diagonal entries
$$(1,1,1,1,e^{2\pi i/3},e^{-2\pi i/3},\ga^2,\overline{\ga}^2,\ga^5,\overline{\ga}^5,\ga^6,\overline{\ga}^6),$$
which are $39$th roots of unity.

\subsubsection{Fusion rules}\label{Haagfusion}
The fusion rules are obtained from the $S$ matrix via the Verlinde
formula.  The fusion matrices for the objects $\pi_1$, $\pi_2$ and
$\mu_1$ are:

{\tiny

$$N_{\pi_1}= \left( \begin {array}{cccccccccccc}
0&1&0&0&0&0&0&0&0&0&0&0
\\\noalign{\medskip}1&1&1&1&1&1&1&1&1&1&1&1\\\noalign{\medskip}0&1&2&1
&1&1&1&1&1&1&1&1\\\noalign{\medskip}0&1&1&2&1&1&1&1&1&1&1&1
\\\noalign{\medskip}0&1&1&1&2&1&1&1&1&1&1&1\\\noalign{\medskip}0&1&1&1
&1&2&1&1&1&1&1&1\\\noalign{\medskip}0&1&1&1&1&1&0&1&1&1&1&1
\\\noalign{\medskip}0&1&1&1&1&1&1&0&1&1&1&1\\\noalign{\medskip}0&1&1&1
&1&1&1&1&0&1&1&1\\\noalign{\medskip}0&1&1&1&1&1&1&1&1&0&1&1
\\\noalign{\medskip}0&1&1&1&1&1&1&1&1&1&0&1\\\noalign{\medskip}0&1&1&1
&1&1&1&1&1&1&1&0\end {array} \right) $$

$$N_{\pi_2}=\left( \begin {array}{cccccccccccc}
0&0&1&0&0&0&0&0&0&0&0&0
\\\noalign{\medskip}0&1&2&1&1&1&1&1&1&1&1&1\\\noalign{\medskip}1&2&2&1
&1&1&1&1&1&1&1&1\\\noalign{\medskip}0&1&1&1&2&2&1&1&1&1&1&1
\\\noalign{\medskip}0&1&1&2&1&2&1&1&1&1&1&1\\\noalign{\medskip}0&1&1&2
&2&1&1&1&1&1&1&1\\\noalign{\medskip}0&1&1&1&1&1&1&1&1&1&1&1
\\\noalign{\medskip}0&1&1&1&1&1&1&1&1&1&1&1\\\noalign{\medskip}0&1&1&1
&1&1&1&1&1&1&1&1\\\noalign{\medskip}0&1&1&1&1&1&1&1&1&1&1&1
\\\noalign{\medskip}0&1&1&1&1&1&1&1&1&1&1&1\\\noalign{\medskip}0&1&1&1
&1&1&1&1&1&1&1&1\end {array} \right) $$

}

and

{\tiny
$$
N_{\mu_1}= \left(\begin{array}{cccccccccccc} 0&0&0&0&0&0&0&1&0&0&0&0\\\noalign{\medskip}
 0&1&1&1&1&1&1&0&1&1&1&1\\\noalign{\medskip}
  0&1&1&1&1&1&1&1&1&1&1&1\\\noalign{\medskip}
0&1&1&1&1&1&1&1&1&1&1&1\\\noalign{\medskip}
0&1&1&1&1&1&1&1&1&1&1&1\\\noalign{\medskip}
0&1&1&1&1&1&1&1&1&1&1&1\\\noalign{\medskip}
0&1&1&1&1&1&1&1&1&1&0&0\\\noalign{\medskip}
1&0&1&1&1&1&1&1&0&1&1&1\\\noalign{\medskip}
0&1&1&1&1&1&1&0&1&0&1&1\\\noalign{\medskip}
0&1&1&1&1&1&1&1&0&1&1&0\\\noalign{\medskip}
0&1&1&1&1&1&0&1&1&1&0&1\\\noalign{\medskip}
0&1&1&1&1&1&0&1&1&0&1&1
\end{array}\right)$$

}

The remaining fusion matrices can be obtained from these three by
permuting the rows and columns.  The specific permutation that
effects these similarity transformations are deduced by comparing
the rows of the $S$ matrix above.  For example, $N_{\sigma_i}$ is
obtained from $N_{\pi_2}$ by the transposition
$\pi_2\leftrightarrow\sigma_i$, since this transposition
transforms the row of the $S$ matrix labeled by $\pi_2$ with the
row labeled by $\sigma_i$.
 Similarly, $N_{\mu_2}$ can be obtained from $N_{\mu_1}$ by a permutation
of rows and columns, specifically, $\mu_1\leftrightarrow\mu_2$,
$\mu_3\leftrightarrow\mu_4$ and $\mu_5\leftrightarrow\mu_6$
converts $N_{\mu_1}$ to $N_{\mu_2}$. The required permutation is
not always order $2$, for example, $N_{\mu_3}$ is obtained from
$N_{\mu_1}$ via the permutation (written in cycle notation)
$(\mu_1,\mu_3,\mu_5)(\mu_2,\mu_4,\mu_6)$.

From these fusion rules we obtain the following:
\begin{prop}\label{haagfusionprop}
The modular category $\mZ(\CH)$ has no non-trivial tensor subcategories.
\end{prop}
\begin{proof}
First observe that $\pi_2\ot\pi_2$ contains every simple object.
Since the diagonal entries of $N_{\pi_2}$ are all positive except
for the trivial object, we conclude that $\pi_2$ appears in $X\ot
X$ for \emph{every} non-trivial simple object in $\mZ(\CH)$ since
every object is self-dual. Thus the tensor subcategory generated
by \emph{any} non-trivial simple object is all of $\mZ(\CH)$.
\end{proof}

\section{MTCs from quantum groups}
\label{s:appl}

From any simple Lie algebra $\g$ and $q\in\C$ with $q^2$ a
primitive $\ell$th root of unity one may construct a ribbon
category $\CC(\g,q,\ell)$ (see \emph{e.g.} \cite{BK}).  One can
also construct such categories from semisimple $\g$, but the
resulting category is easily seen to be a direct product of those
constructed from simple $\g$.  We shall say these categories (or
direct products of them) are of \textbf{quantum group type}. There
is an (often overlooked) subtlety concerning the degree $\ell$ of
$q^2$ and the unitarizability of the category $\CC(\g,q,\ell)$.
Let $m$ denote the maximal number of edges between any two nodes
of the Dynkin diagram for $\g$ with $\g$ simple, so that $m=1$ for
Lie types $ADE$, $m=2$ for Lie types $BCF_4$ and $m=3$ for Lie
type $G_2$.  Provided $m\mid \ell$, $\CC(\g,q,\ell)$ is a unitary
category for $q=e^{\pi i/\ell}$ (see \cite{wenzlcstar}).  If
$m\nmid\ell$, this is not always true  and in fact there is
usually no choice of $q$ to make $\CC(\g,q,\ell)$ unitary (see
\cite{rowellmathz} and \cite{grunit}).  In \cite{finkel} it is
shown that the fusion category associated with level $k$
representations of the affine Kac-Moody algebra $\hat{\g}$ is
tensor equivalent to $\CC(\g,q,\ell)$ for $\ell=m(k+h_\g)$ where
$h_\g$ is the dual Coxeter number.  In these cases the categories
are often denoted $(X_r,k)$, where $\g$ is of Lie type $X$ with
rank $r$ and $k=\ell/m-h_\g$ is the level. We will use this
abbreviated notation except when $m\nmid\ell$.

Our
goal is to prove the following:

\begin{theorem}\label{noquantum}
The modular categories $\mZ(\CE)$ and $\mZ(\CH)$ are not monoidally equivalent to any
category of quantum group type.
\end{theorem}

Before we proceed to the proof, we give a few more details on the
categories $\CC(\g,q,\ell)$ with $\g$ simple.  The (isomorphism
classes of) objects in $\Cg$ are labeled by a certain finite
subset $C_\ell$ of the dominant weights of $\g$.  The size
$|C_\ell|$ is the \emph{rank} of $\Cg$; for $\g$ simple,
generating functions for $|C_\ell|$ are found in \cite{rowell06}.
For any object $X$, we have that $X\cong X^*$ if and only if the
corresponding simple $\g$-module $V$ satisfies $V\cong V^*$, in
which case we say that $X$ is \emph{self-dual}.  If every object
is self-dual, we will say the category itself is self-dual.  A
simple object $X_\la$ is self-dual if and only if $-w_0(\la)=\la$,
where $w_0$ is the longest element of the Weyl group of $\g$, (see
e.g. \cite{GoodWall}[Exercise 5.1.8.4]).  In a ribbon category
(such as $\Cg$) it is always true that $X^{**}\cong X$ for any
object $X$, so non-self dual objects always appear in pairs.
Observe also that the unit object $\1$ is always self-dual.  The
twists for the simple objects in $\Cg$ are powers of $q^{1/M}$
where $M$ is the order of the quotient group $P/Q$ of the weight
lattice $P$ by the root lattice $Q$.  In particular, $M=(r+1)$ for
$A_r$, and for all other Lie types $M\leq 4$.

\begin{proof}
Observe that $\mZ(\CE)$ has rank $10$ and has exactly one pair of
simple non-self-dual objects, and $8$ simple self-dual objects (up
to isomorphism). For $\mZ(\CE)$ the statement follows from the
following fact (to be established below): \emph{no rank $10$
category of quantum group type has exactly one pair of simple
non-self-dual objects.}

It is immediate that $\mZ(\CE)$ cannot be a direct product
$\F\boxtimes\T$ of two non-trivial modular categories.  First
notice that one may assume that $\rank(\F)=2$ and $\rank(\T)=5$.
Since rank $2$ modular category are all self-dual and simple
non-self-dual objects appear in pairs, $\F\boxtimes\T$ has either
$0$, $4$ or $8$ simple non-self-dual objects, while $\mZ(\CE)$ has
exactly $2$.

Next we observe that if $\g$ is of Lie type $A_1$, $B_r$, $C_r$,
$D_{2t}$, $E_7$, $E_8$, $F_4$ or $G_2$, all of the simple objects
in the corresponding category are self-dual, since $-1$ is the
longest element of the Weyl group.  So we may immediately
eliminate categories of these Lie types from consideration. This
leaves only Lie types $A_r$ ($r\geq 2$), $D_r$ ($r\geq 5$ and odd)
and $E_6$ as possibilities.

From \cite{rowell06} we have the following generating functions
for $|C_\ell|$, the rank of $\Cg$:
\begin{enumerate}
\item $A_r$: $\frac{1}{(1-x)^{r+1}}=\sum_{k=0}^\infty
\binom{r+k}{k}x^k$. \item $D_r$:
$\frac{1}{(1-x)^4(1-x^2)^{r-3}}=1+4x+(r+7)x^2+(8+4r)x^3+\cdots$
\item $E_6$:
$\frac{1}{(1-x)^3(1-x^2)^3(1-x^3)}=1+3x+9x^2+20x^3+42x^4+\cdots$
\end{enumerate}
where the coefficient of $x^k$ is the rank of $\CC(\g,q,\ell)$,
where $\ell=k+h$ with $h$ the dual Coxeter number of the root
system of $\g$.

To determine the rank $10$ type $A$ categories, we must solve
$\binom{r+k}{k}=10$ for $(r,k)$.  The only positive integer
solutions are $(1,9)$, $(9,1)$, $(3,2)$ and $(2,3)$.  These
correspond to $(A_1,9)$ and $(A_9,1)$ (i.e. at $11$th roots of
unity) and
 $(A_2,3)$ and $(A_3,2)$ (i.e. at $6$th roots of unity).  The category
$(A_1,9)$ has only self-dual objects, and the remaining three categories
each have at least $4$ non-self-dual objects.  For type $D_r$
with $r\geq 5$ and $E_6$ we see that no rank $10$ categories
appear.

Next let us consider the (self-dual) rank 12 category $\mZ(\CH)$.  Since $\mZ(\CH)$ has no tensor
subcategories by Prop. \ref{haagfusionprop}, $\mZ(\CH)$ is not the product of two modular categories.
Using the generating
functions from \cite{rowell06} we determine all rank $12$ self-dual quantum group categories.  The
following pairs $(X_r,k)$ are the rank $r$ Lie type $X$ quantum groups
at level $k$ that have exactly $12$ simple objects:

$$\{(G_2,5),(A_1,11),(B_8,2),(C_{11},1),(D_5,2),(E_7,3)\}.$$

In addition the categories $\CC(\mathfrak{so}_5,q,9)$ and
$\CC(\mathfrak{so}_{11},q,13)$ have rank $12$. The twists
$\theta_i$ in $\mZ(\CH)$ include $13$th roots of unity, so that we
may immediately eliminate all of the above except $(A_1,11),
(C_{11},1)$ and $\CC(\mathfrak{so}_{11},q,13)$.  By a level-rank
duality theorem in \cite{rowellmathz} the pairs $(A_1,11)$ and
$\CC(\mathfrak{so}_{11},q,13)$ have the same fusion rules, and
moreover each contains a tensor subcategory eliminating them from
consideration.  To eliminate $(C_{11},1)$ we must work a little
harder.  The category $\CC(\mathfrak{sp}_{22},e^{\pi i/26},26)$
contains a simple object $X:=X_{(1,0,\cdots,0)}$ corresponding to
the $22$ dimensional representation of $\mathfrak{sp}_{22}$. The
quantum dimension of this object is
$$\frac{[11][24]}{[1][12]}=2\cos(\pi/13)=1.94188\ldots,$$ where
$[n]$ is the usual $q$-number at $q=e^{\pi i/26}$.  Since the
simple objects in $\mZ(\CH)$ have quantum dimensions in
$\{1,3d,3d+1,3d+2\},$ where $d=\frac{3+\sqrt{13}}{2}>3$, we see
that $\mZ(\CH)$ cannot be obtained from $(C_{11},1)$. Thus
$\mZ(\CH)$ is not a quantum group type category.
\end{proof}

\section{Doubled Categories}\label{double}

Another way in which modular categories may be constructed is as
the quantum double $\mZ(\CC)$ of a spherical fusion category
\cite{MugerII}.  In this section we will prove the following:

\begin{theorem}\label{nodouble}
The modular categories $\mZ(\CE)$ and $\mZ(\CH)$ are not braided
monoidally equivalent to the double of any braided fusion
category.
\end{theorem}

\begin{proof}

First observe that  by \cite{Kas}[Corollary XIII.4.4] and
\cite{MugerII}[Lemma 7.1] the double $\mZ(\CC)$ of a braided
fusion category $\CC$ contains a braided tensor subcategory
equivalent to $\CC$.  In the case of $\mZ(\CH)$, we showed in
Prop. \ref{haagfusionprop} that the only tensor subcategories are
the trivial subcategory and $\mZ(\CH)$ itself. So $\mZ(\CH)$ is
not the double of \emph{any} braided fusion category.

Suppose that $\mZ(\CE)$ is braided monoidally equivalent to the
double $\mZ(\CC)$ of some braided fusion category $\CC$, then
$\CC$ is equivalent to some braided fusion subcategory
$\CC^\prime\subset \mZ(\CE)$. Since $\mZ(\CE)$ is modular we may
further assume that $\CC^\prime$ is pre-modular and that
$\dim(\CC^\prime)^2=\dim(\mZ(\CE))=(6+2\sqrt{3})^2$ by
\cite{MugerII}[Theorem 1.2].  Thus $\dim(\CC^\prime)=6+2\sqrt{3}$
which implies that there exists some simple object
$X\in\CC^\prime$ with $\dim(X)\neq 1$.  Let $X$ be such an object,
then $\dim(X)\in\{1+\sqrt{3},2+\sqrt{3},3+\sqrt{3}\}$ since $X$
would be a simple object of $\mZ(\CE)$. The inequality
$\dim(\1)^2+\dim(X)^2\leq\dim(\CC^\prime)=6+2\sqrt{3}$ implies
that $\dim(X)=1+\sqrt{3}$, and this forces $\CC^\prime$ to have 3
simple objects of dimension $1,1$ and $1+\sqrt{3}$.  But it is
known (\cite{Ostrik1}) that no such category can be braided,
therefore, $\mZ(\CE)$ cannot be a double of any braided fusion
category.
\end{proof}

\section{$\SL(2,\Z)$ Image}\label{sl2z}

Every TQFT gives rise to a projective representation of $\SL(2,\Z)$ via
 $$\begin{pmatrix}0 & -1\\1 &
0\end{pmatrix}\rightarrow S, \quad \begin{pmatrix}1 & 1\\0 &
1\end{pmatrix}\rightarrow T.$$  If the TQFT has a corresponding
RCFT, then the resulting representation of $SL(2,\Z)$ has a finite
image group, and the kernel is a congruence subgroup
\cite{Bantay}\cite{Xu3}. It is an open question if this is true
for every TQFT. In particular, any TQFT whose representation of
$SL(2,\Z)$ has an infinite image or has a non-congruence subgroup
kernel is not a CSW TQFT.  But representations of $SL(2,\Z)$ from
$\mZ(\CE)$ and $\mZ(\CH)$ both behave as those of TQFTs from
RCFTs.

\begin{theorem}

\begin{enumerate}

\item The representation of $SL(2,\Z)$ from $\mZ(\CE)$ has a finite
image in $\U(10)$;

\item The representation of $SL(2,\Z)$ from $\mZ(\CH)$ has a finite
image in $\U(12)$, and its kernel is a congruence subgroup.

\end{enumerate}

\end{theorem}

\begin{proof}
First let us consider the category $\mZ(\CE)$.
We wish to show that the $10$-dimensional unitary representation
of $SL(2,\Z)$ given by
$$\begin{pmatrix}0 & -1\\1 &
0\end{pmatrix}\rightarrow S, \quad \begin{pmatrix}1 & 1\\0 &
1\end{pmatrix}\rightarrow T$$ has finite image, where $S$ and $T$
are as in Section \ref{SandT}.  We accomplish this as follows: Let
$G=\lan S,T\ra$ be the group generated by $S$ and $T$.  Observe
that we immediately have the following relations in $G$, since $T$
is a diagonal matrix whose nonzero entries are $12$th roots of
unity:
\begin{eqnarray}\label{strels}
S^4=T^{12}=I,\quad (S T)^3=S^2.
\end{eqnarray}

Additionally, we find the following relation:
\begin{equation}\label{Arel}
(T^4S T^6S)^6=I.
\end{equation}  In fact, $A:=(T^4S T^6S)^2$
is a diagonal order $3$ matrix and will play an important role in
what follows.

Now consider the normal closure, $N$, of the cyclic subgroup
generated by $A$ in $G$. We need to see how $G$ acts on $N$, and
for this it is enough to understand the action of $S$ and $T$ on a
set of generators of $N$.  We will employ the standard notation
for conjugation in a group: $g^h:=hgh^{-1}$.  Defining $B:=A^S$,
$C:=B^T$ and $D:=C^S$, we find that a set of generators for $N$ is
$\{A,B,C,D\}$.  This is established by determining the conjugation
action of $S$ and $T$ on these generators as follows:
\begin{eqnarray}\label{staction}
&A^T=A,&\quad A^S=B,\quad B^T=C,\quad B^S=A,\\
 &C^T=D,&\quad C^S=D,\quad D^T=B,\quad D^S=C.\nonumber
\end{eqnarray}
Thus the subgroup generated by $\{A,B,C,D\}$ is normal and
contains $A$, hence is equal to $N$. Furthermore, one has the
following relations in $N$:
\begin{eqnarray}\label{Nrels}
A^3=I,\quad B^A=D,\quad C^A=B,\quad D^A=C,\quad D^B=A
\end{eqnarray}
which are sufficient to determine all other conjugation relations
among the generators of $N$.  Having established these relations
in $G$, we proceed to analyze the structure of the abstract group
$\widehat{G}$ on generators $\{s,t,a,b,c,d\}$ satisfying relations
(\ref{strels}), (\ref{Arel}), (\ref{staction}) and (\ref{Nrels})
(where $S$ is replaced by $s$ etc.). Clearly $G$ is a quotient of
$\widehat{G}$ since these relations hold in $G$.

First observe $N$ is a quotient of the abstract group
$$\widehat{N}:=\lan a,b,c,d\ |\ a^3=I,b^a=d,c^a=b,d^a=c, d^b=a\ra$$
and that $\widehat{N}\lhd\widehat{G}$. We compute $|\widehat{N}|=24$, so that
if $\widehat{H}:=\widehat{G}/\widehat{N}$ is a finite group, then $\widehat{G}$ is
finite hence $G$ is finite. Next observe that $$\widehat{H}=\lan s,t\
|\ s^4=t^{12}=(t^4st^6s)^2=I,(st)^3=s^2\ra,$$
 \emph{i.e.} with relations (\ref{strels}) together with $(t^4st^6s)^2=I$,
 and we have a short exact sequence
 $$1\rightarrow\widehat{N}\rightarrow\widehat{G}\rightarrow\widehat{H}\rightarrow 1.$$
Using MAPLE, we find that $|\widehat{H}|=1296$ so that
$|\widehat{G}|=(24)(1296)=31104$.  Thus $|G|$ is a finite group of
order dividing $31104=(2^7)(3^5)$.  Note that since $|G|$ is not divisible by
$10$ we can conclude that the representation above is reducible.

Incidentally, the group $\SL(2,\Z/36\Z)$ (i.e. $\SL(2,\Z)$ with
entries taken modulo $36$) has order $31104$ just as $\widehat{G}$
above does.  While we expect that the kernel of $\SL(2,\Z)\rightarrow G$
is a congruence subgroup, it is \emph{not} true that
$\widehat{G}\cong\SL(2,\Z/36\Z)$, according to computations with GAP (\cite{GAP}).
However, both of these groups are solvable, and their (normal) Sylow-$2$ and Sylow-$3$ subgroups
\emph{are} isomorphic.  This implies that $\widehat{G}$ and $\SL(2,\Z/36\Z)$
are semi-direct products of the same two groups.  This of course does not imply that the kernel
of $\SL(2,\Z)\rightarrow G$ is non-congruence, merely that the obvious guess is incorrect.

For $\mZ(\CH)$ we use the $12\times 12$ matrices $S$ and $T$ found
in Section \ref{HaagSandT}. For $N$ odd, the following generators
and relations for $\PSL(2,\Z/N\Z)$ are found in \cite{Sunday}:
$$(AB)^3=A^2=B^N=(B^4AB^{\frac{N+1}{2}}A)^2=I.$$
Setting $A=S$, $B=T$ and $N=39$, we easily verify these relations, so that the image of $SL(2,\Z)$
for the Haagerup MTC $\mZ(\CH)$ is finite and
the kernel is a congruence subgroup.

\end{proof}

\section{Representation of the braid groups}\label{braidreps}

Representations of the braid group $\B_n$ can be obtained from any simple object $X$
in a braided tensor category $\CC$.  The construction is as follows.  For $1\leq i\leq n-1$, Let $\beta_i$ be the usual
generators of the braid group satisfying
$\beta_i\beta_{i+1}\beta_i=\beta_{i+1}\beta_i\beta_{i+1}$ and $\beta_i\beta_j=\beta_j\beta_i$ for $|i-j|>1$.
The braiding operator
$c_{X,X}\in\End(X^{\ot 2})$ acts on $\End(X^{\ot 2})$ by composition, and we define invertible operators
in $\End(X^{\ot n})$ by:
$\phi_X^n(\beta_i)=Id_X^{\ot i-1}\ot c_{X,X}\ot Id_X^{\ot n-i-1}$.
This defines a representation $\B_n\rightarrow \GL(\End(X^{\ot n}))$ by
$\beta_i.f\rightarrow \phi_X^n(\beta_i)\circ f$ which is unitarizable in case
$\CC$ is a unitary ribbon category.

These representations are rarely irreducible.  In fact, if $Y$ is a
simple subobject of $X^{\ot n}$, then $\Hom(Y,X^{\ot n})$ is
(isomorphic to) a $\B_n$-subrepresentation of $\End(X^{\ot n})$,
since $\Hom(Y,X^{\ot n})$ is obviously stable under composition
with $\phi_X^n(\beta_i)$. However, the $\B_n$-subrepresentations
of the form $\Hom(Y,X^{\ot n})$ are not irreducible for all $Y$
unless the algebra $\End(X^{\ot n})$ is generated by
$\{\phi_X^n(\beta_i)\}$.  It is a technically difficult problem to
determine the irreducible constituents of $\End(X^{\ot n})$ as a
$\B_n$-representation; few general techniques are available. One
useful criterion is the following proposition \cite{TuWe}[Lemma
5.5]:
\begin{prop}\label{TW5.5}
Suppose $X$ is a simple self-dual object in a ribbon category $\CC$, such that
\begin{enumerate}
\item[(a)] $X^{\ot 2}$ decomposes as a direct sum of $d$ distinct simple objects $X_i$ and
\item[(b)] $\phi_X^2(\beta_1)$ has $d$ distinct eigenvalues.
\end{enumerate}
Then $\B_3$ acts irreducibly on $\Hom(X,X^{\ot 3})$.
\end{prop}
A generalization of this result to spaces of the form
$\Hom(Y,X^{\ot 3})$ and for repeated eigenvalues would be of
considerable value.

An important question for a given MTC $\CC$ is the following: do these representations
of $\B_n$ factor over finite groups for all $X$ and all $n$,
or is there a choice of $X$ so that the image of $\B_n$ is infinite (say, for all $n\geq 3$)?

\begin{prop}

Representations of the braid groups $\B_n, n\geq 3$ from the
simple object $X_4$ of $\mZ(\CE)$ has a dense image in the
projective unitary group.

\end{prop}

\begin{proof}
For the simple object $X_4$ of $\mZ(\CE)$, we have $X_4^{\ot
2}=\1+X_4+V$.  So $\Hom(\mathbf{1},X_4^2)$, $\Hom(X_4,X_4^2)$,
$\Hom(V,X_4^2)$ are all $1$-dimensional.  Set $q=e^{\pi i/6}$. The
eigenvalues of braiding $\phi_{X_4}^2(\beta_i)$ are computed as:
$a=\frac{1-\sqrt{3}i}{2}=q^{-2},b=\frac{-\sqrt{3}+i}{2}=-q^{-1},c=\frac{\sqrt{3}+i}{2}=q$.
Since these eigenvalues are distinct, it follows from Prop.
\ref{TW5.5} that the representation is irreducible. The projective
order of braid generators of $\B_3$ is $12$ because $12$ is the
smallest $m$ so that $a^m=b^m=c^m$.  It follows from
\cite{LRW}[Prop. 6.8] that the image of the representation of
$\B_n$ for $n\geq 3$ afforded by $\mZ(\CE)$ with each braid strand colored by
$X_4$ is infinite, and dense in the projective unitary group.
\end{proof}

For the Haagerup MTC $\mZ(\CH)$ we cannot conclude the image is infinite without considerably
more work.  The necessary techniques are somewhat ad hoc and go beyond the scope of this paper.
We plan to give an account of these techniques in a subsequent
article.

We give a brief explanation of the difficulties one encounters in this case. The smallest non-trivial
representation of $\B_3$ is $7$ dimensional, for example acting on the vector space
$\Hom(\mu_2,\mu_1^{\ot 3})$.  Set $\ga=e^{2\pi i/13}$ as above. The eigenvalues of the
operators $(\phi_{\mu_1}^3(\beta_i))^2$ restricted to this space are
$$\ga^4\cdot\{1,1,e^{\pm 2\pi i/3},\ga^{\pm 2},\ga^{-5}\}$$
which can easily be computed from the twists.
However, we do not know that this representation is irreducible, or indeed, not a
sum of 1-dimensional representations, since Prop. \ref{TW5.5} only applies to spaces
of the form $\Hom(X,X^{\ot 3})$.  A $10$ dimensional representation of $\B_3$ with
the right form is $\Hom(\mu_1,\mu_1^{\ot 3})$,
and the corresponding eigenvalues of $(\phi^3_{\mu_1}(\beta_i))^2$ are:
$$\ga^4\{1,1,1,e^{\pm 2\pi i/3},\ga^{\pm 2},\ga^{-5},\ga^{\pm 6}\}.$$
But the eigenvalues of $\phi^3_{\mu_1}(\beta_i)$ are some choices of
square roots of these values which will clearly not be distinct.

The technique for showing that the image is infinite is as
follows.  First find an irreducible subrepresentation of dimension
$d$.  Next verify that the corresponding image is not imprimitive
by checking the ``no-cycle condition" of \cite{LRW} or by some other means.  Then check
that the projective order of the images of $\beta_i$ does not
occur for primitive linear groups of degree $d$ by checking the
lists in \cite{Feit1} and \cite{Feit2} of primitive linear groups of degree $d\leq 10$. For example, if $\B_3$ acts primitively
on some $d$-dimensional
subrepresentation $W$ of $\Hom(\mu_1,\mu_1^{\ot 3})$ and has a) $2\leq d$
and b) both $3$rd and $13$th roots of unity occur as eigenvalues of $(\phi^3_{\mu_1}(\beta_i))^2$ acting on $W$ then the image must
be infinite.

\section{Central Charge and Orbifold/Coset
Constructions}\label{coset}

Since we do not know how to cover general compact Lie groups $G$
in the quantum group setting, we will restrict our discussion to
the semi-simple cases.  The orbifold and coset CFTs for WZW models
with semi-simple Lie groups $G$ have been constructed
mathematically.  Although complete analysis of all possible
orbifold and coset candidates for $\mZ(\CE)$ and $\mZ(\CH)$ seems
impossible, we will give evidence that orbifold and coset
constructions are unlikely to realize them.

If a TQFT comes from a RCFT, then a relation between
the topological central charge of the TQFT and the chiral central
charge of the corresponding RCFT exists.  The topological
central charge of a TQFT is defined as follows:

\begin{defn}

Let $d_i$ be the quantum dimensions of all simple types $X_i,
i=1,\cdots, n$ of an MTC $\CC$, $\theta_i$ be the twists, define
the \textit{total quantum order} of $\CC$ to be
$D=\sqrt{\sum_{i=1}^n d_i^2}$, and $D_{+}=\sum_{i=1}^n \theta_i
d_i^2$.  Then $\frac{D_{+}}{D}=e^{\frac{\pi i c}{4}}$ for some
rational number $c$.  The rational number $c$ defined modulo $8$
will be called the topological central charge of the MTC $\CC$.

\end{defn}

Each CSW TQFT corresponds to a RCFT.  The chiral central charge of
the RCFT is a rational number $c_v$.  We have the following:

\begin{prop}

If a TQFT has a corresponding RCFT, then $c_v=c$ mod $8$, in
particular this is true for CSW TQFTs.

\end{prop}

\begin{proof}
This relation first appeared in \cite{rehren}.
For another explanation, see
\cite{Kitaev2} on Page 66.  For general
unitary TQFTs, it is not known if the boundary theories are always
RCFTs, and it is an open question if there is a similar identity.
See the references in \cite{Kitaev2} and \cite{Witten2}.

\end{proof}

Since $\mZ(\CE), \mZ(\CH)$ have topological central charge $=0$, a
corresponding CFT, if exist, would have chiral central charge $=0$
mod $8$. To rule out the possibility of coset and orbifold
constructions of $\mZ(\CE)$ or $\mZ(\CH)$, we need to have a list
of all chiral central charge $=0$ mod $8$ CFTs, and their
orbifolds. This question seems hard. So we will only consider, as
examples, the case of CFTs with chiral central charge $0$ or $24$.
Even with this restriction, the problem is still hard, so we will
further restrict our discussion to simple quantum group
categories, i.e. those from simple Lie algebras plus their
orbifolds and certain cosets.

As shown in \cite{MugerO}, the orbifold construction in CFT cannot
be formulated purely in a categorical way (coset construction has
not been attempted systematically in the categorical framework).
In the case of quantum group categories, this problem can be
circumvented by the following detour: the corresponding CFTs are
WZW models, coset and orbifold CFTs of WZW models are
mathematically constructed (see \cite{MugerO}\cite{Xu2} and the
references therein). We will then take the corresponding MTCs of
the resulting CFTs as the cosets or orbifolds of the quantum group
categories.

We collect some facts about orbifold and coset CFTs that we need
in this section from \cite{MugerO}\cite{Xu2}\cite{CFT}. Given a
simple Lie algebra $\mathfrak{g}$ and a level $k$, the WZW CFT has
chiral central charge $c=\frac{k\cdot \textrm{dim}
\mathfrak{g}}{k+h_{\mathfrak{g}}}$, where $h_{\mathfrak{g}}$ is
the dual Coxter number of $\mathfrak{g}$.  Given a CFT with a
discrete finite automorphism group $G$ on its chiral algebra $A$,
then the orbifold CFT based on $\textrm{Rep}(A^G)$ of the fixed
algebra has the same chiral central charge.

The coset construction is complicated and is defined for any pair
of Lie groups $H\subset G$. We will restrict ourselves to the
cases that $\mathfrak{p}\subset \mathfrak{g}$ such that both are
simple, and $\mathfrak{p}$ is an isolated maximal subalgebra as in
Tables $2$ and $5$ of \cite{Bais}.  The total quantum order $D_{G/H}$ of
a coset $H\subset G$ MTC is $D_{G/H}=\frac{D_G}{D_H} \cdot
d^2(G/H)$, where $d^2(G/H)$ is the index of type $II_1$ subfactors
\cite{Xu2}.  By Jones' celebrated theorem \cite{Jones}, if
$d(G/H)\leq 2$, then $d(G/H)=2cos(\pi/r)$ for some $r\geq 3$.

Given a simple Lie algebra $\mathfrak{g}$, and a level $k$. Let
$\mathfrak{p}$ be a simple subalgebra, and $\chi$ be the Dynkin
embedding index of $\mathfrak{p}$ in $\mathfrak{g}$. Then the
central charge of the resulting CFT $\mathfrak{g}/\mathfrak{p}$ is
$c_{\mathfrak{g}/\mathfrak{p}}=\frac{k\cdot
\textrm{dim}\mathfrak{g} }{k+h_{\mathfrak{g}}} - \frac{\chi \cdot
k \cdot \textrm{dim} \mathfrak{p}}{\chi \cdot k
+h_{\mathfrak{p}}}.$ Recall the Dynkin embedding index for a pair
of simple Lie algebras $\mathfrak{p}\subset \mathfrak{g}$: let
$\lambda$ be a highest weight of $\mathfrak{g}$, then
$\lambda=\oplus_{\mu\in P_{+}} b_{\lambda\mu} \mu$, where $P_{+}$
is the set of dominant weights of $\mathfrak{p}$, then
$\chi_{\mathfrak{g}/\mathfrak{p}}=\sum_{\mu\in P_{+}}
b_{\lambda\mu} \frac{\chi_{\mu}}{\chi_{\lambda}}$, where
$\chi_{\lambda}=\frac{\textrm{dim} \lambda\cdot
(\lambda,\lambda+2\rho_{\mathfrak{g}})}{2\cdot
\textrm{dim}\mathfrak{g}}, \chi_{\mu}=\frac{\textrm{dim} \mu\cdot
(\mu,\mu+2\rho_{\mathfrak{p}})}{2\cdot \textrm{dim}\mathfrak{p}}$.

\begin{prop}

\begin{enumerate}

\item If the total quantum order of a unitary MTC $\CC$ from a CFT is
$D$, then any nontrivial orbifold of $\CC$ has total quantum order $\geq
2D$.

\item A chiral central charge $24$ unitary CFT from a simple Lie
algebra is one of the following:
$$(A_6,7), (A_{24},1), (B_{12},2), (C_4,10), (D_{24},1).$$
Moreover, neither $\mZ(\CE)$ nor $\mZ(\CH)$ is an orbifold of
those CFTs.

\item The only chiral central charge $24$ unitary coset CFT of the
form $\mathfrak{g}/\mathfrak{p}$ for simple
$\mathfrak{p},\mathfrak{g}$ in Tables $2$ and $5$ of \cite{Bais}
is from the embedding $A_7\subset D_{35}$ with embedding index
$\chi=10$, and $D_{35}$ is at level $k=2$.  The resulting coset
TQFT is neither $\mZ(\CE)$ nor $\mZ(\CH)$.

\end{enumerate}

\end{prop}

\begin{proof}

(1):  Let $D_A$ be the total quantum order of the MTC corresponds to CFT
$A$, then the orbifold MTC has total quantum order $|G|\cdot D_A$, and
the inequality follows.

(2): A chiral unitary central charge $0$ CFT is trivial, and the
orbifolds of the trivial CFT are (twisted) quantum double of
finite groups whose quantum dimensions are all integers \cite{BK}.
But we know $\mZ(\CE)$ and $\mZ(\CH)$ both have non-integral
quantum dimensions, hence they are not orbifolds of the trivial
CFT.

In \cite{Sc}, $71$ CFTs of chiral central charge $24$ are listed.
A simple inspection gives our list for simple algebras.  More
directly, we can find the list by solving Diophantine equations
$24=\frac{k\cdot \textrm{dim} \mathfrak{g}}{k+h_{\mathfrak{g}}}$
for all simple Lie algebras.

$(A_{24},1)$ corresponds to $\SU(25)$ at a $26$th root of unity.
This is a rank $25$ abelian theory, with all categorical
dimensions of simple objects equal to $1$.  Similarly,
corresponding to $(D_{24},1)$ is an abelian rank $4$ category.  So
any orbifold theory will have global quantum dimension $N^2\cdot 25$ or
$N^2\cdot 4$ for some integer $N$ which is obviously not
$(6+2\sqrt{3})^2$ or $(\frac{39+9\sqrt{13}}{2})^2$.

$(B_{12},2)$ corresponds to $SO(25)$ at a $50$th root of unity.
Since the global quantum dimension of this rank $16$ category must reside
in $\Q[e^{\pi i/100}]$, it is clear that no integer multiple of
its global quantum dimension can be $(6+2\sqrt{3})^2$ or
$(\frac{39+9\sqrt{13}}{2})^2$.

$(C_4,10)$ corresponds to $Sp(8)$ at a $30$th root of unity,
having rank $1001$ \cite{rowell06}.  Since the quantum dimension
of any simple object is $\geq 1$ for a unitary theory, thus the
total quantum order of any orbifold theory is at least $2\sqrt{1001}\cong
63.3$. Similarly, $(A_6,7)$ is a unitary MTC from $\SU(7)$ at a
$14$th root of unity.  By \cite{rowell06}, its rank is ${13
\choose 7}$.  It follows that the total quantum order of $(A_6,7)\geq
\sqrt{{13 \choose 7}}$. Hence any nontrivial orbifold of $(A_6,7)$
will have a total quantum order $D_G \geq 2 \sqrt{{13 \choose 7}}\cong
82.8$. But $D_{\mZ(\CE)}=6+2\sqrt{3}$, and
$D_{\mZ(\CH)}=\frac{39+9\sqrt{13}}{2}\cong 35.7$, hence it is
impossible for either to be an orbifold of $(C_4,10)$ or
$(A_6,7)$.

(3):  The coset TQFT is obtained from $(D_{35})_2/(A_7)_{20}$.
This embedding is as follows:  the fundamental representation
$\mu=\omega_4$ of $SU(8)$ is of dimension $70$.  $\omega_4$ has a
symmetric invariant bilinear form which gives rise to the
embedding of $SU(8)$ into $SO(70)$, corresponding to the
fundamental representation $\lambda=\omega_1$ of $D_{35}$.  Hence
the branching rule for $\lambda$ is simply $\mu$, and the coset
theory has a simple object labeled by $(\lambda,\mu)$.  The
embedding index can be computed using the formula above:
$\chi_{\mu}=\frac{70\cdot 18}{2\cdot 63}$,
$\chi_{\lambda}=\frac{70\cdot 69}{2\cdot 35\cdot 69}$, so
$\chi_{\lambda/\mu}=10$. For level $k=1$ of $D_{35}$, this
embedding is conformal, i.e. the resulting coset has chiral
central charge $0$. For level $k=2$, the resulting coset has
chiral central charge $24$. By the formulas $(18.42)$ on Page
$805$ \cite{CFT} (cf. \cite{Xu2}), the twist of the simple object
$(\lambda,\mu)$ in the coset is
$\frac{\theta_{\lambda}}{\theta_{\mu}}$.  When $k=2$, $D_{35}$
corresponds to $SO(70)$ at a $70$th root of unity $q=e^{\frac{\pi
i}{70}}$, and the twist of $\lambda$ is
$\theta_{\lambda}=q^{(\lambda,\lambda+2\rho)}=e^{-\frac{\pi
i}{70}}$. The level for $A_7$ is $20$ since the embedding index is
$10$. So$(A_7)_{20}$ corresponds to $SU(8)$ at a  $28$th root of
unity $a=e^{\frac{\pi i}{28}}$, so the twist of $\mu$ is
$\theta_{\mu}=a^{(\mu,\mu+2\rho)}=e^{\frac{9\pi i}{14}}$. The
twists of $\mZ(\CE)$ are all $12$th root of unity, and $\mZ(\CH)$
all $39$th root of unity. Since the ratio
$\theta_{(\lambda,\mu)}=\frac{\theta_{\lambda}}{\theta_{\mu}}=e^{-\frac{23
\pi i}{35}}$ can never be a $12$th or $39$th root of unity, hence
this coset MTC is neither $\mZ(\CE)$ nor $\mZ(\CH)$.

\end{proof}

\section{Appendix}

\subsection{Category $\frac{1}{2} E_6$}\label{original cat}

The category $\CE=\frac{1}{2}E_6$ is a unitary monoidal spherical
category of rank $3$. The following is the information for its
structure. (Details can be found in \cite{HH})

$\bullet$ simple objects:

$\{1,x,y\}$

$\bullet$ fusion rule:

$x^2=1+2x+y, xy=x=yx$

$\bullet$ basis:

$v^{1}_{11} \in V^{1}_{11}$, $v^{x}_{1x} \in V^{x}_{1x}$,
$v^{x}_{x1} \in V^{x}_{x1}$,$v^{y}_{1y} \in V^{y}_{1y}$,
$v^{y}_{y1} \in V^{y}_{y1}$, $v^{x}_{xy} \in V^{x}_{xy}$,
$v^{x}_{yx} \in V^{x}_{yx}$, $v^{1}_{yy} \in V^{1}_{yy}$,
$v^{1}_{xx} \in V^{1}_{xx}$, $v^{y}_{xx} \in V^{y}_{xx}$, and
$v_1$, $v_2 \in V^{x}_{xx}$, where $V^z_{xy}$ denotes
$\Hom_{\CE}(xy,z)$.

$\bullet$ associativities:

$a^{y}_{y,y,y}=$ $a^{x}_{x,y,y}=$ $a^{x}_{y,y,x}=$
$a^{1}_{x,y,x}=$ $a^{1}_{x,x,y}=$ $a^{y}_{x,x,y}=$
$a^{1}_{y,x,x}=$ $a^{y}_{y,x,x}=1$, $a^{y}_{x,y,x}=$
$a^{x}_{y,x,y}= -1$, $a^{x}_{x,y,x}=
\left[\begin{array}{cc}1&0\\0&-1\end{array} \right]$,
$a^{x}_{x,x,y}= \left[\begin{array}{cc}0&i\\-i&0\end{array}
\right]$, $a^{x}_{y,x,x}=
\left[\begin{array}{cc}0&1\\1&0\end{array} \right]$,
$a^{1}_{x,x,x}= \frac{1}{\sqrt{2}} e^{7 \pi i
 /12} \left[\begin{array}{cc}1&i\\1&-i\end{array}
\right]$, $a^{y}_{x,x,x}= \frac{1}{\sqrt{2}} e^{7 \pi i
 /12} \left[\begin{array}{cc}i&1\\-i&1\end{array}
\right]$

$a^{x}_{x,x,x}=\left[ \begin{smallmatrix}
\frac{-1+\sqrt{3}}{2}&\frac{-1+\sqrt{3}}{2}& \frac{1- \sqrt{3}}{4}
e^{\pi i /6}&\frac{1- \sqrt{3}}{4} e^{2\pi i /3}
&\frac{1- \sqrt{3}}{4} e^{2\pi i /3}&\frac{1- \sqrt{3}}{4} e^{\pi i /6} \\
\frac{-1+\sqrt{3}}{2}&\frac{1-\sqrt{3}}{2}& \frac{1- \sqrt{3}}{4}
e^{\pi i /6}&\frac{1- \sqrt{3}}{4} e^{2\pi i /3}&-\frac{1-
\sqrt{3}}{4} e^{2\pi i /3}&-\frac{1- \sqrt{3}}{4} e^{\pi i /6}
\\1&1&-\frac{1}{2} (e^{\pi i /6} -1)&\frac{1}{2} e^{5\pi i /6}&\frac{1}{2} (e^{-\pi i /3} + i)&\frac{1}{2} e^{\pi i /3}
\\1&1&\frac{1}{2} e^{\pi i /3}&\frac{1}{2} (e^{-\pi i /3} + i)&\frac{1}{2} e^{5\pi i /6}&-\frac{1}{2} (e^{\pi i /6} -1)
\\1&-1&-\frac{1}{2} (e^{\pi i /6} -1)&\frac{1}{2} e^{5\pi i /6}&-\frac{1}{2} (e^{-\pi i /3} +
i)&-\frac{1}{2} e^{\pi i /3}
\\-1&1&-\frac{1}{2} e^{\pi i /3}&-\frac{1}{2} (e^{-\pi i /3} + i)&\frac{1}{2} e^{5\pi i /6}&-\frac{1}{2} (e^{\pi i /6} -1)
\end{smallmatrix}
\right]$

$\bullet$ notations for the dual basis:

we use $v^{xy}_{z} \in \Hom_{\CE}(z,xy)$ to denote the dual basis
of $v^{z}_{xy}$ in the sense that $v^{z}_{xy} \circ v^{xy}_{z} =
id_z$, and use $v^{1}$ and $v^{2}$ for dual bases of $v_1$ and
$v_2$, respectively.

$\bullet$ rigidity:

$d_y:=v^1_{yy}$, $b_y:=v^{yy}_1$, $d_x:=v^1_{xx}$,
$b_x:=(1+\sqrt{3})v^{xx}_1$.

$\bullet$ quantum dimensions:

$\dim_{\CE}(1)=1$, $\dim_{\CE}(y)=1$, $\dim_{\CE}(x)=1+\sqrt{3}$.

\subsection{Definitions and Lemmas}

In this section, we follow Section $3$ of \cite{MugerII}.

\begin{defn}
\label{h.b.def} Let $C$ be a strict monoidal category and let $x
\in C$. A half braiding $e_x$ for $x$ is a family $\{e_x(y) \in
\Hom_C(xy,yx),y \in C \}$ of isomorphisms satisfying

(i) Naturality: $f \otimes id_x \circ e_x(y) =e_x(z) \circ id_x
\otimes f \:\:\:\forall f:y \rightarrow z.$

(ii) The braid relation: $e_x(y\otimes z)=id_y \otimes e_x(z)
\circ e_x(y) \otimes id_z \:\:\:\forall y,z \in C$.

(iii) Unit property: $e_x(1)=id_x$.
\end{defn}

The following lemma is equivalent to Lemma 3.3 of \cite{MugerII}.

\begin{lemma}\label{only consider simple}
Let $C$ be semisimple and $\{ x_i,i \in \Gamma \}$ a basis of
simple objects. Let $z \in C$. Then there is a one-to-one
correspondence between (i) families of morphisms $\{e_z(x_i) \in
\Hom_C(zx_i,x_iz),i \in \Gamma \}$ such that

$e_z(x_k) \circ id_z \otimes f = f \otimes id_z \circ id_{x_i}
\otimes e_z(x_j) \circ e_z(x_i) \otimes id_{x_j}\:\:\: \forall
i,j,k \in \Gamma, f \in \Hom_C(x_ix_j,x_k)$,

and (ii) families of morphisms $\{ e_z(x) \in \Hom_C(zx,xz,x \in C
\}$ satisfying 1. and 2. from the Definition~\ref{h.b.def}. All
$e_z(x),x\in C$ are isomorphisms iff all $e_z(x_i),i \in \Gamma$
are isomorphisms.
\end{lemma}

\begin{defn}
\label{qd} The quantum double $\mZ(C)$ of a strict monoidal category
$C$ has as objects pairs $(x,e_x)$, where $x \in C$ and $e_x$ is a
half braiding. The morphisms are given by

$\Hom_{\mZ(C)}((x,e_x),(y,e_y))= \{ f \in \Hom_C(x,y)| id_z \otimes
f\circ e_x(z)=e_y(z) \circ f \otimes id_z \:\:\: \forall z \in C
\}$.

The tensor product of objects is given by $(x,e_x) \otimes (y,e_y)
= (xy,e_{xy})$, where

$e_{xy}(z)=e_x(z) \otimes id_y \circ id_x \otimes e_y(z)$.

The tensor unit is $(1,e_1)$ where $e_1(x)=id_x$. The composition
and tensor product of morphisms are inherited from $C$. The
braiding is given by

$c((x,e_x),(y,e_y))=e_x(y)$.
\end{defn}

\subsection{Solutions for the half braiding}

For any object $x,y \in \CE$, $\Hom_{\CE}(xy,yx)$ has a basis
consisting of morphisms of the type $(v^{yx}_{x_k})_j \circ
(v^{x_k}_{xy})_i$ where $k \in \Gamma$ and $1 \leq i, j \leq
\dim(\Hom_{\CE}(xy,x_k))$. We parameterize each half braiding
as a linear combination of such basis vectors and need to
determine the coefficients satisfying all constraints in
Definition~\ref{h.b.def}. However, from Lemma~\ref{only consider
simple}, we only need to consider naturality with respect to the
basis morphisms in Section~\ref{original cat}. The following are
the solutions where $x \in \CE$ has $5$ half braidings denoted by
$e_{x_i}, i=1,2,\cdots,5$:

$--\;\;$$e_y(y)=-v^{yy}_{1}\circ v^{1}_{yy}$

$--\;\;$$e_y(x)=i v^{xy}_{x} \circ v^{x}_{yx}$

$--\;\;$$e_{x_1}(y)=i v^{yx}_{x} \circ v^{x}_{xy}$

$--\;\;$$e_{x_1}(x)=iv^{xx}_{1} \circ v^{1}_{xx}+ v^{xx}_{y}\circ
v^{y}_{xx}+ e^{-\pi i/3}v^1\circ v_1  +e^{-5\pi i/6} v^2\circ v_2$

$--\;\;$$e_{x_2}(y)=i v^{yx}_{x} \circ v^{x}_{xy}$

$--\;\;$$e_{x_2}(x)=e^{-5\pi i/6}v^{xx}_{1} \circ
v^{1}_{xx}+e^{2\pi i/3} v^{xx}_{y}\circ
v^{y}_{xx}+\frac{1-\sqrt{3}}{2} v^1\circ
v_1+\left(\frac{\sqrt{3}}{2}\right)^{1/2}i v^2\circ v_1+
\left(\frac{\sqrt{3}}{2}\right)^{1/2}v^1\circ
v_2+\frac{\sqrt{3}-1}{2}i v^2\circ v_2$

$--\;\;$$e_{x_3}(y)=i v^{yx}_{x} \circ v^{x}_{xy}$

$--\;\;$$e_{x_3}(x)=e^{-5\pi i/6}v^{xx}_{1} \circ
v^{1}_{xx}+e^{2\pi i/3} v^{xx}_{y}\circ
v^{y}_{xx}+\frac{1-\sqrt{3}}{2} v^1\circ
v_1-\left(\frac{\sqrt{3}}{2}\right)^{1/2}i v^2\circ v_1-
\left(\frac{\sqrt{3}}{2}\right)^{1/2}v^1\circ
v_2+\frac{\sqrt{3}-1}{2}i v^2\circ v_2$

$--\;\;$$e_{x_4}(y)=-i v^{yx}_{x} \circ v^{x}_{xy}$

$--\;\;$$e_{x_4}(x)=e^{-\pi i/3}v^{xx}_{1} \circ v^{1}_{xx}+
e^{\pi i/6}v^{xx}_{y}\circ v^{y}_{xx}+ \frac{1}{\sqrt{2}}e^{\pi
i/4}v^1\circ v_1+\frac{1}{\sqrt{2}}e^{-\pi i/4} v^2\circ v_1+
\frac{1}{\sqrt{2}}e^{\pi i/4}v^1\circ v_2+
\frac{1}{\sqrt{2}}e^{3\pi i/4}v^2\circ v_2$

$--\;\;$$e_{x_5}(y)=-i v^{yx}_{x} \circ v^{x}_{xy}$

$--\;\;$$e_{x_5}(x)=e^{2\pi i/3}v^{xx}_{1} \circ v^{1}_{xx}+
e^{-5\pi i/6}v^{xx}_{y}\circ v^{y}_{xx}+
\frac{1}{\sqrt{2}}e^{-3\pi i/4}v^1\circ
v_1+\frac{1}{\sqrt{2}}e^{-\pi i/4} v^2\circ v_1+
\frac{1}{\sqrt{2}}e^{\pi i/4}v^1\circ v_2+
\frac{1}{\sqrt{2}}e^{-\pi i/4}v^2\circ v_2$

$--\;\;$$e_{1+x}(y)=v^{y1}_{y}\circ v^{y}_{1y}-i v^{yx}_{x}\circ
v^{x}_{xy}$

$--\;\;$$e_{1+x}(x)=(-2+\sqrt{3})v^{x1}_{x} \circ
v^{x}_{1x}+(2\sqrt{3}-3)v^{1}\circ
v^{x}_{1x}+(2\sqrt{3}-3)v^{2}\circ v^{x}_{1x}+e^{-5 \pi
i/6}v^{x1}_{x}\circ v_1 +e^{-\pi i/3}v^{x1}_{x}\circ v_2 +
v^{xx}_{1}\circ v^{1}_{xx}+i v^{xx}_{y}\circ
v^{y}_{xx}+\frac{\sqrt{3}-1}{\sqrt{2}}e^{-5\pi i/12}v^1\circ v_1
+\frac{\sqrt{3}-1}{\sqrt{2}}e^{3\pi i/4}v^2\circ v_1
+\frac{\sqrt{3}-1}{\sqrt{2}}e^{-3\pi i/4}v^1\circ v_2
+\frac{\sqrt{3}-1}{\sqrt{2}}e^{\pi i/12}v^2\circ v_2 $

$--\;\;$$e_{y+x}(y)=-v^{yy}_{1}\circ v^{1}_{yy}-i v^{yx}_{x}\circ
v^{x}_{xy}$

$--\;\;$$e_{y+x}(x)=(-2+\sqrt{3})iv^{xy}_{x}\circ
v^{x}_{yx}+(2\sqrt{3}-3)iv^{1}\circ
v^{x}_{yx}-(2\sqrt{3}-3)iv^{2}\circ v^{x}_{yx}+e^{\pi
i/6}v^{xy}_{x}\circ v_{1}+e^{-\pi i/3}v^{xy}_{x}\circ
v_{2}-v^{xx}_{1}\circ v^{1}_{xx}-iv^{xx}_{y}\circ
v^{y}_{xx}+\frac{\sqrt{3}-1}{\sqrt{2}}e^{7\pi i/12}v^1\circ v_1
+\frac{\sqrt{3}-1}{\sqrt{2}}e^{3\pi i/4}v^2\circ v_1
+\frac{\sqrt{3}-1}{\sqrt{2}}e^{-3\pi i/4}v^1\circ v_2
+\frac{\sqrt{3}-1}{\sqrt{2}}e^{-11\pi i/12}v^2\circ v_2 $

$--\;\;$$e_{1+y+x}(y)=-v^{y1}_{y}\circ v^{y}_{1y}
+v^{yy}_{1}\circ v^{1}_{yy} +iv^{yx}_{x}\circ v^{x}_{xy}$

$--\;\;$$e_{1+y+x}(x)=2e^{-5\pi i/6}v^{xy}_{x}\circ
v^{x}_{1x}+2e^{-5\pi i/6}v^{1}\circ v^{x}_{1x}+2e^{\pi
i/6}v^{2}\circ v^{x}_{1x}+\frac{2-\sqrt{3}}{2}e^{5\pi
i/6}v^{x1}_{x}\circ v^{x}_{yx}+\frac{\sqrt{3}-1}{\sqrt{2}}e^{-\pi
i/4}v^{1}\circ v^{x}_{yx}+\frac{\sqrt{3}-1}{\sqrt{2}}e^{-\pi
i/4}v^{2}\circ v^{x}_{yx}+v^{xx}_{1}\circ
v^{1}_{xx}-iv^{xx}_{y}\circ
v^{y}_{xx}+\frac{\sqrt{3}-1}{4}iv^{x1}_{x}\circ
v_{1}+\frac{\sqrt{3}-1}{4}v^{x1}_{x}\circ
v_{2}+\frac{1}{\sqrt{2}}e^{-7\pi i/12}v^{xy}_{x}\circ
v_{1}+\frac{1}{\sqrt{2}}e^{-\pi i/12}v^{xy}_{x}\circ
v_{2}+\frac{\sqrt{3}-1}{\sqrt{2}}e^{5\pi i/12}v^{1}\circ
v_{1}+\frac{\sqrt{3}-1}{\sqrt{2}}e^{-\pi i/12}v^{2}\circ v_{2}$

We will use the following notation:

$\mathbf{1}:=(1,e_1)$, $Y:=(y,e_y)$, $X_i:=(x,e_{x_i})$ for
$i=1,2,\cdots,5$, $U:=(1+x,e_{1+x})$, $V:=(y+x,e_{y+x})$, and
$W:=(1+y+x,e_{1+y+x})$.

It is not hard to see that all these $10$ objects are simple in
the quantum doubled category $\mZ(\CE)$, and not isomorphic to each
other by considering each $\Hom_{\mZ(\CE)}$-space in
Definition~\ref{qd}. Furthermore, these $10$ objects completes the
list of representatives of isomorphism classes of simple objects
in $\mZ(\CE)$ by the fact $\dim \mZ(\CE)=(\dim \CE)^2$ (see Theorem 4.14
of \cite{MugerII}).

To decompose each tensor product into direct sum of simple
objects, we need to compute fusion morphisms satisfying the
conditions in Definition~\ref{qd}. After parameterizing each
morphism as a linear combination of basis morphisms in
~\ref{original cat}, to find solutions for each coefficient is
purely algebraic computation, from which we can determine the
dimension of each $\Hom_{\mZ(\CE)}$-space.  This can be done
easily.

\end{document}